\documentclass[12pt]{article}

\def\ra{\rightarrow}

\def\ss{\subseteq}

\def\P{{\cal P}}
\def\L{{\cal L}}

\def\zbar{\overline{z}}
\def\zetabar{\overline{\zeta}}
\def\O{\Omega}
\def\Sz{Szeg\H{o}}
\def\P{{\cal P}}
\def\B{{\cal B}}

\def\L{{\cal L}}

\def\dbar{\overline{\partial}}

 \def\HollowBox #1#2{{\dimen0=#1 \advance\dimen0 by -#2       
       \dimen1=#1 \advance\dimen1 by #2                       
        \vrule height #1 depth #2 width #2                    
        \vrule height 0pt depth #2 width #1                   
        \llap{\vrule height #1 depth -\dimen0 width \dimen1}%
       \hskip -#2                                             
       \vrule height #1 depth #2 width #2}}                   
 \def\BoxOpTwo{\mathord{\HollowBox{6pt}{.4pt}}\;}             

\def\endpf{\hfill $\BoxOpTwo$}

\font\teneufm=eufm10
\font\seveneufm=eufm7
\font\fiveeufm=eufm5
\newfam\eufmfam
\textfont\eufmfam=\teneufm
\scriptfont\eufmfam=\seveneufm
\scriptscriptfont\eufmfam=\fiveeufm

\newfam\msbfam
\font\tenmsb=msbm10  scaled \magstep1 \textfont\msbfam=\tenmsb
\font\sevenmsb=msbm7 scaled \magstep1 \scriptfont\msbfam=\sevenmsb
\font\fivemsb=msbm5  scaled \magstep1 \scriptscriptfont\msbfam=\fivemsb
\def\Bbb{\fam\msbfam \tenmsb}

\def\CC{{\Bbb C}}

\newtheorem{theorem}{Theorem}[section]

\newtheorem{proposition}[theorem]{Proposition}
\newtheorem{lemma}[theorem]{Lemma}
\newtheorem{remark}[theorem]{Remark}

\newtheorem{definition}{Definition}[section]

\usepackage{graphicx}

\usepackage{amsmath}

\begin{document}

\begin{center}
{\Large \bf On a Construction of L. Hua for Positive
\bigskip \\
Reproducing Kernels}
\bigskip \bigskip \\
Steven G. Krantz\footnote{Supported in part by a grant from
the National Science Foundation and a grant from the Dean of
Graduate Studies at Washington University.}
\end{center}
\vspace*{.25in}

\begin{quote}
{\bf Abstract:}  \sl
We study a positive reproducing kernel for holomorphic functions
on complex domains.  This kernel, which induces what has now come
to be known as the Berezin transform, is manufactured from the Bergman
kernel using an idea of L. K. Hua.  The kernel has important
analytic and geometric properties which we develop in some detail.
\end{quote}

\section{Introduction}

Let $\Omega \ss \CC^n$ be a bounded domain (i.e., a connected open set).
Following the general rubric of ``Hilbert space with reproducing kernel''
laid down by Nachman Aronszajn [ARO], both the
Bergman space $A^2(\Omega)$ and the Hardy space $H^2(\Omega)$ have
reproducing kernels.  We shall provide the details of these assertions
below.

The Bergman kernel (for $A^2$) and the Szeg\H{o} kernel (for $H^2$) both have
the advantage of being canonical.  But neither is positive, and this
makes them tricky to handle.  The Bergman kernel can be treated with the theory
of the Hilbert integral (see [PHS]) and the Szeg\H{o} kernel can often
be handled with a suitable theory of singular integrals (see [KRA2]).

It is a classical construction of Hua (see [HUA]) that one can
use the Szeg\H{o} kernel to produce another reproducing kernel
${\cal P}(z, \zeta)$ which also reproduces $H^2$ but which is
positive. In this sense it is more like the Poisson kernel of
harmonic function theory. In point of fact, this so-called
Poisson-Szeg\H{o} kernel coincides with the Poisson kernel
when the domain is the disc $D$ in the complex plane $\CC$.
Furthermore, the Poisson-Szeg\H{o} kernel solves the Dirichlet
problem for the invariant Laplacian (i.e., the
Laplace-Beltrami operator for the Bergman metric) on the ball
in $\CC^n$. Unfortunately a similar statement about the
Poisson-Szeg\H{o} kernel cannot be made on any other domain
(although we shall explore substitute results on strongly
pseudoconvex domains in the present paper).

We want to develop these ideas with the Szeg\"{o} kernel
replaced by the Bergman kernel. This notion was developed
independently by Berezin [BER] in the context of quantization
of K\"{a}hler manifolds. Indeed, one assigns to a bounded
function on the manifold the corresponding Toeplitz operator.
This process of a assigning a linear operator to a function is
called {\it quantization}. A nice exposition of the ideas
appears in [PET]. Further basic properties may be found in
[ZHU].

Approaches to the Berezin transform are often
operator-theoretic (see [ENG1], [ENG2]), or sometimes
geometric [PET]. Our point of view here will be more
function-theoretic. We shall repeat (in perhaps new language)
some results that are known in other contexts. And we shall
also enunciate and prove new results. We hope that the mix
serves to be both informative and useful.

It is a pleasure to thank M. Engli\v{s} and R. Rochberg for helpful
conversations.

\section{Fundamental Ideas}

If $\O \ss \CC^n$ is a bounded domain then set
$$
A^2(\O) = \left \{f \ \hbox{holomorphic on} \ \O: \int_\O |f(z|^2 \, dV(z) < \infty \right \} \, .
$$
[Here $dV$ is standard Euclidean volume measure on $\O$.]
Of course $A^2$ is equipped with the inner product
$$
\langle f, g \rangle_{A^2(\O)} = \int_\O f(z) \overline{g(z)} \, dV(z) \, .
$$
Then $A^2$ is a subspace of $L^2(\O)$, and it can be shown (see [KRA3]) that
$A^2(\O)$ is a Hilbert space.  We have the fundamental lemma:

\begin{lemma} \sl
Let $K \ss \O$ be a compact subset of $\O \ss \CC^n$.  There is a constant $C = C(K, n)$
such that, if $f \in A^2(\O)$, then
$$
\sup_{z \in K} |f(z)| \leq C \cdot \|f\|_{A^2(\Omega)} \, .
$$
\end{lemma}

We shall not prove the lemma here, but refer the reader to [KRA3] or [KRA4]
for the details.

Now if $z \in \O$ is a fixed point then, applying the lemma
with $K = \{z\}$, we find that the linear functional
$$
e_z: A^2(\O) \ni f \longmapsto f(z)
$$
is bounded.  By the Riesz representation theorem, there then exists
a function $k_z \in A^2(\O)$ such that
$$
f(z) = e_z(f) = \langle f, k_z \rangle
$$
for all $f \in A^2(\O)$.  We set $K(z,\zeta) = K_\O(z, \zeta) = \overline{k_z(\zeta)}$ and
write
$$
f(z) = \int_\Omega K(z,\zeta) f(\zeta) \, dV(\zeta) \, .
$$
This is the Bergman reproducing formula and $K(z, \zeta)$ is the
Bergman (reproducing) kernel.  

There is a similar theory for $H^2$.  Fix a bounded domain $\O$.  
Define
$$
H^2(\O) = \{f \ \hbox{holomorphic on} \ \O: |f|^2 \ \hbox{has a harmonic majorant on} \ \O\} \, .
$$
This definition is equivalent to several other natural definitions of $H^2$; see [KRA1] for
the details.  In particular, it can be shown that an $H^2$ function $f$ has an $L^2(\partial \Omega)$
boundary function $\widetilde{f}$ and that $f$ is the Poisson integral of $\widetilde{f}$.
It is convenient to set $\|f\|_{H^2(\Omega)} = \|\widetilde{f}\|_{L^2(\partial \Omega)}$.  This
definition of the norm is equivalent to several other standard definitions---again see [KRA1].

We now have the fundamental lemma:

\begin{lemma} \sl
Let $K \ss \O$ be a compact subset of $\O \ss \CC^n$.  There is a constant $C' = C'(K, n)$
such that, if $f \in H^2(\O)$, then
$$
\sup_{z \in K} |f(z)| \leq C \cdot \|f\|_{H^2(\Omega)} \, .
$$
\end{lemma}
Again, details of the proof are omitted.

As a consequence, if a point $z \in \O$ is fixed, then we can be sure
that the functional
$$
e'_z: H^2(\O) \ni f \longmapsto f(z)
$$     
is bounded.  By the Riesz representation theorem, there then exists
a function $k'_z \in A^2(\O)$ such that
$$
f(z) = e'_z(f) = \langle f, k'_z \rangle
$$
for all $f \in H^2(\O)$.  We set $S(z,\zeta) = S_\O(z, \zeta) = \overline{k'_z(\zeta)}$ and
write
$$
f(z) = \int_{\partial \Omega} S(z,\zeta) f(\zeta) \, d\sigma(\zeta) \, .
$$
[Here $d\sigma$ is standard area measure (i.e., {\it Hausdorff measure}) on $\partial \Omega$.]
This is the \Sz\ reproducing formula and $S(z, \zeta)$ is the \Sz\ (reproducing) kernel. 

Of course the projection
$$
P_B: L^2(\Omega) \ra A^2(\O)
$$
is well defined by
$$
P_B f(z) = \int_\Omega K(z, \zeta) f(\zeta) \, dV(\zeta) \, .
$$
Likewise the projection
$$
P_S: L^2(\partial \Omega) \ra H^2(\O)
$$
is well defined by
$$
P_S f(z) = \int_{\partial \O} S(z, \zeta) f(\zeta) \, d\sigma(\zeta) \, .
$$
These two facts establish the centrality and importance of the kernels
$K$ and $S$.  But neither kernel is positive, and that makes their
analysis difficult.

\section{Positive Kernels}

In the seminal work [HUA], L. Hua proposed a program for producing a positive
kernel from a canonical kernel.  He defined
$$ 
\P(z, \zeta) = \frac{|S(z, \zeta)|^2}{S(z,z)} \, ,
$$
where $S$ is the standard \Sz\ kernel on a given bounded domain $\O$.  Now we have

\begin{proposition} \sl
Let $\O$ be a bounded domain with $C^2$ boundary and $S$ its \Sz\ kernel.
With $\P(z, \zeta)$ as defined above, and with $f \in C(\overline{\O})$ holomorphic
on $\Omega$, we have
$$
f(z) = \int_{\partial \Omega} \P(z, \zeta) f(\zeta) \, d\sigma(\zeta)
$$
for all $z \in \Omega$.
\end{proposition}
{\bf Proof:}  Fix $z \in \Omega$.  Define $g(\zeta) = \overline{S(z,\zeta)} \cdot f(\zeta)/S(z,z)$.
Then it is easy to see that $g \in H^2(\Omega)$ as a function of $\zeta$.  As a result,
\begin{align}
\int_{\partial \Omega} f(\zeta) {\cal P}(z, \zeta) \, d\sigma(\zeta) & = 
         \int_{\partial \Omega} \left [ f(\zeta) \cdot \frac{\overline{S(z,\zeta)}}{S(z,z)} \right ] \cdot S(z, \zeta) \, d\sigma(\zeta) \notag \\
             & =  \int_{\partial \Omega} g(\zeta) \cdot S(z, \zeta) \, d\sigma(\zeta) \notag \\
	     & =  g(z)  \notag \\
	     & =  f(z) \, .  \tag*{$\BoxOpTwo$}	  \\ \notag
\end{align}

Notice that the fact that $f$ is continuous on $\overline{\O}$ is used to
guarantee that $g \in H^2$. It is natural to ask whether the result of the
proposition extends to all functions $f \in H^2(\O)$. For this, it would
suffice to show that $C(\overline{\O}) \cap {\cal O}(\O)$ is dense in
$H^2(\O)$. In fact this density result is known to be true when $\O$ is
either strongly pseudoconvex or of finite type in the sense of
Catlin/D'Angelo/Kohn. In fact one can reason as follows (and we thank
Harold Boas for this argument): Let $f \in H^2(\Omega)$. Then certainly $f
\in L^2(\partial \Omega)$ and, just by measure theory, one can approximate
$f$ in $L^2$ norm by a function $\varphi \in C^\infty(\partial \Omega)$.
Let $\Phi = P_S \varphi$, the \Sz\ projection of $\varphi$. Then, since $P_S$
is a continuous operator on $L^2(\partial \Omega)$, the function $\Phi$ is
an $L^2(\partial \Omega)$ approximant of $f$. But it is also the case, by regularity theory
of the $\overline{\partial}_b$ operator, that $\Phi = P_S\varphi$ is in
$C^\infty(\overline{\O})$. That proves the needed approximation result. Of
course a similar argument would apply on any domain on which the \Sz\
projection maps smooth functions to smooth functions. See [STE] for some
observations about this matter.

Now Hua did not consider his construction for the Bergman kernel, but in
fact it is just as valid in that context.  We may define
$$
{\cal B}(z, \zeta) = \frac{|K(z, \zeta)|^2}{K(z,z)} \, .
$$
We call this the {\it Poisson-Bergman kernel}.  Then we have

\begin{proposition} \sl
Let $\O$ be a bounded domain and $K$ its Bergman kernel.
With ${\cal B}(z, \zeta)$ as defined above, and with $f \in C(\overline{\O})$ holomorphic
on $\Omega$, we have
$$
f(z) = \int_{\partial \Omega} {\cal B}(z, \zeta) f(\zeta) \, dV(\zeta)
$$
for all $z \in \Omega$.
\end{proposition}
The proof is just the same as that for Proposition 3.1, and we omit the details.
One of the purposes of the present paper is to study properties of
the Poisson-Bergman kernel $\B$.

Of course the Poisson-Bergman kernel is real, so it will also reproduce the real
parts of holomorphic functions.  Thus, in one complex variable, the integral
reproduces harmonic functions.  In several complex variables, it reproduces pluriharmonic
functions.

Again, it is natural to ask under what circumstances Proposition 3.2 holds
for all functions in the Bergman space $A^2(\Omega)$.  The question
is virtually equivalent to asking when the elements that are continuous
on $\overline{\Omega}$ are dense in $A^2$.  Catlin [CAT] has given
an affirmative answer to this query on any smoothly bounded pseudoconvex
domain.

One of the features that makes the Bergman kernel both important and useful
is its invariance under bihlomorphic mappings.  This fact is useful in conformal
mapping theory, and it also gives rise to the Bergman metric.  The fundamental result is
this:

\begin{proposition} \sl
Let $\Omega_1, \Omega_2$ be domains in $\CC^n.$
Let $f: \Omega_1 \ra \Omega_2$ be biholomorphic.  Then
$$
\mbox{\rm det}\, J_\CC f(z)  K_{\Omega_2}(f(z),f(\zeta)) \mbox{\rm det}\, \overline{J_\CC f(\zeta)}
                         = K_{\Omega_1}(z,\zeta) . 
$$
\end{proposition}
Here $J_\CC f$ is the complex Jacobian matrix of the mapping $f$.  Refer to [KRA1] and
[KRA4] for more on this topic.

It is useful to know that the Poisson-Bergman kernel satisfies a similar transformation law:

\begin{proposition}
Let $\Omega_1, \Omega_2$ be domains in $\CC^n.$
Let $f: \Omega_1 \ra \Omega_2$ be biholomorphic.  Then
$$
\B_{\Omega_2}(f(z),f(\zeta)) |\mbox{\rm det}\, J_\CC f(\zeta)|^2
                         = \B_{\Omega_1}(z,\zeta) \, . 
$$
\end{proposition}
{\bf Proof:}  Of course we use the result of Proposition 3.3.  Now
\begin{eqnarray*}
\B_{\Omega_1}(z, \zeta) & = & \frac{|K_{\Omega_1}(z, \zeta)|^2}{K_{\Omega_1}(z,z)} \\
			& = & \frac{|\hbox{det}\, J_\CC f(z) \cdot K_{\Omega_2}(f(z), f(\zeta)) \cdot \overline{\hbox{det}\, J_\CC f(\zeta)}|^2}{
				 \hbox{det}\, J_\CC f(z) \cdot K_{\Omega_2}(f(z), f(z)) \cdot \overline{\hbox{det}\, J_\CC f(z)}} \\
			& = & \frac{|\hbox{det} \, J_\CC f(\zeta)|^2 \cdot |K_{\Omega_2}(z, \zeta)|^2}{K_{\Omega_2}(f(z), f(z))} \\
			& = & |\hbox{det} \, J_\CC f(\zeta)|^2 \cdot \B_{\Omega_2}(z,z) \, .
\end{eqnarray*}

We conclude this section with an interesting observation about the Berezin transform---see [ZHU].

\begin{proposition} \sl
The operator
$$
\B f(z) = \int_B \B(z, \zeta) f(\zeta) \, dV(\zeta) \, ,
$$
acting on $L^1(B)$, is univalent.
\end{proposition}
{\bf Proof:}   In fact it is useful to take advantage of the symmetry
of the ball.  We can rewrite the Poisson-Bergman integral as
$$
\int_B f \circ \Phi_z (\zeta) \, dV(\zeta) \, ,
$$
where $\Phi_z$ is a suitable automorphism of the ball.  Then it is clear
that this integral can be identically zero in $z$ only if $f \equiv 0$.
That completes the proof.

Another, slightly more abstract, way to look at this matter is
as follows (we thank Richard Rochberg for this idea, and see
also [ENG1]). Let $f$ be any $L^1$ function on $B$.   For $w \in B$
define
$$
g_w(\zeta) = \frac{1}{(1 - \overline{w} \cdot \zeta)^{n+1}} \, .
$$
If $f$ is bounded on the ball, let
$$
T_f: g \mapsto P_B(fg) \, .
$$
We may write the Berezin transform now as
$$
\Lambda f (w,z) = \frac{\langle T_f g_z, g_w \rangle}{\langle g_w, g_w \rangle } \, .
$$
This function is holomorphic in $z$ and conjugate holomorphic in $w$.  The statement
that the Berezin transform $\B f() \equiv 0$ is the same as $\Lambda f(z,z) = 0$.
But it is a standard fact (see [KRA1]) that we may then conclude that $\Lambda f(w,z) \equiv 0$.
But then $T_f g_z \equiv 0$ and so $f \equiv 0$.  So the Berezin transform is univalent.
\endpf
\smallskip \\

\section{Boundary Behavior}

It is natural to want information about the boundary limits of potentials
of the form $\B f$ for $f \in L^2(\O)$.  We begin with a simple lemma:

\begin{lemma} \sl
Let $\O$ be a bounded domain and $\B$ its Poisson-Bergman kernel.  If $z \in \O$ is
fixed, then
$$
\int_\O \B(z, \zeta) \, dV(\zeta) = 1 \, .
$$
\end{lemma}
{\bf Proof:}  Certainly the function $f(\zeta) \equiv 1$ is an element of
the Bergman space on $\O$.  As a result,
$$
1 = f(z) = \int_\O \B(z, \zeta) f(\zeta) \, dV(\zeta) = \int_\O \B(z, \zeta) \, dV(\zeta)
$$
for any $z \in \Omega$.
\endpf
\smallskip \\

Our first result is as follows:

\begin{proposition} \sl
Let $\Omega$ be the ball $B$ in $\CC^n$.
Then the mapping
$$
f \mapsto \int_\O \B(z, \zeta) f(\zeta) \, dV(\zeta) 
$$
sends $L^p(\Omega)$ to $L^p(\Omega)$, $1 \leq p \leq \infty$.
\end{proposition}
{\bf Proof:}  We know from the lemma that
$$
\|\B(z, \ \cdot \ ) \|_{L^1(\Omega)} = 1
$$
for each fixed $z$.  An even easier estimate shows that
$$
\|\B( \ \cdot \ , \zeta ) \|_{L^1(\Omega)} \leq 1
$$
for each fixed $\zeta$.  Now Schur's lemma, or the generalized Minkowski
inequality, give the desired conclusion.
\endpf
\smallskip \\

\begin{proposition} \sl
Let $\O \ss \CC^n$ be the unit ball $B$.
Let $f \in C(\overline{\O})$.  Let $F = \B f$.  Then $F$ extends
to a function that is continuous on $\overline{\O}$.  Moreover, if
$P \in \partial \O$, then
$$
\lim_{\Omega \ni z \ra P} F(z) = f(P) \, .
$$
\end{proposition}
{\bf Proof:}  Let $\epsilon > 0$.  Choose $\delta > 0$ such that
if $z, w \in \overline{\Omega}$ and $|z - w| < \delta$ then
$|f(z) - f(w)| < \epsilon$.  Let $M = \sup_{\zeta \in \overline{\O}} |f(\zeta)|$.
Now, for $z \in \O$, $P \in \partial \O$, 
and $|z - P| < \epsilon$, we have that
\begin{eqnarray*}
|F(z) - f(P)| & = & \left | \int_\O \B(z, \zeta) f(\zeta) \, dV(\zeta) - f(P) \right | \\
	      & = & \left | \int_\O \B(z, \zeta) f(\zeta) \, dV(\zeta) -
	                    \int_\O \B(z, \zeta) f(P) \, dV(\zeta) \right | \\
	   & \leq & \int_{\zeta \in \O \atop
	                     |\zeta - P| < \delta} \B(z, \zeta) |f(\zeta) - f(P)| \, dV(\zeta)  \\
              && \quad + \int_{\zeta \in \O \atop
		             |\zeta - P| \geq \delta} \B(z, \zeta) |f(\zeta) - f(P)| \, dV(\zeta) \\
	   & \leq & \int_{\zeta \in \O \atop
	                     |\zeta - P| < \delta} \B(z, \zeta) \cdot \epsilon \, dV(\zeta) +
		    \int_{\zeta \in \O \atop
		             |\zeta - P| \geq \delta} \B(z, \zeta) \cdot 2M \, dV(\zeta) \\
	   & \equiv & I + II \, .
\end{eqnarray*}

Now the lemma tells us that $I = \epsilon$.   Also we know that the Poisson-Bergman kernel
for the ball is
$$
\B(z, \zeta) = c_n \frac{(1 - |z|^2)^{n+1}}{|1 - z \cdot \overline{\zeta}|^{2n+2}} \, .
$$
Thus, by inspection, $\B(z, \zeta) \ra 0$ as $z \ra P$ for $|\zeta - P| \geq \delta$.
Thus $II$ is smaller than $\epsilon$ as soon as $z$ is close enough to $P$.

In summary, for $z$ sufficiently close to $P$, $|F(z) - f(P)| < 2\epsilon$.  That is what
we wished to prove.
\endpf
\smallskip \\

Arazy and Engli\v{s} have in fact shown that the last result is true
on any pseudoconvex domain for which each boundary point is a peak
point (for the algebra $A(\Omega)$ of functions continuous on the closure and
holomorphic inside).  Thus the result is true in particular on strongly
pseudoconvex domains (see [KRA1]) and finite type domains in $\CC^2$ (see [BEF]).

Here is another way to look at the matter on strongly pseudoconvex domains.  In fact our observation, at
the end of the proof of the last proposition, about the vanishing of $\B(z, \zeta)$ for
$z \ra P$ and $|\zeta - P| \geq \delta$ is a tricky point and not generally known.  On
a strongly pseudoconvex domain $\O$ we have Fefferman's asymptotic expansion [FEF].  This
says that, in suitable local holomorphic coordinates near a boundary point $P$, we have
$$
K_\O(z, \zeta) = \frac{c_n}{(1 - z \cdot \overline{\zeta})^{n+1}} +
                 k(z,\zeta) \cdot \log  |1 - z \cdot \overline{\zeta}| \, .  \eqno (4.3)
$$
Thus, using an argument quite similar to the one that we carry out in detail in Section 5 for
the Poisson-\Sz\  kernel, one can obtain an asymptotic expansion for the Poisson-Bergman kernel.
One sees that, in local coordinates near the boundary.
$$
\B_\O(z, \zeta) = c_n \cdot \frac{(1 - |z|^2)^{n+1}}{|1 - z \cdot \overline{\zeta}|^{2n+2}} + {\cal E}(z, \zeta) \, ,
$$
where ${\cal E}$ is a kernel that induces a smoothing operator.  In particular, the singularity of ${\cal E}$ will
be measurably less than the singularity of the lead term.  So it will still be the case
that $\B(z,\zeta) \ra 0$ as $z \ra P \in \partial \Omega$ and $|\zeta - P| \geq \delta$.  So we have:

\begin{proposition} \sl
Let $\O \ss \CC^n$ be a smoothly bounded, strongly pseudoconvex domain in $\CC^n$.  Let
$f \in C(\overline{\O})$.  Then the function $\B f$ extends to be continuous on $\overline{\O}$.
Moreover, if $P \in \partial \Omega$, then
$$
\lim_{\O \ni z \ra P} \B f(z) = f(P) \, .
$$
\end{proposition}

It is natural, from the point of view of measure theory and harmonic analysis, to want
to extend the result of Proposition 4.3 to a broader class of functions.   To this
end we introduce a maximal function to use as a tool.

\begin{definition} \rm
Let $\Omega$ be a smoothly bounded, strongly pseudoconvex domain in $\CC^n$.  If $z, \zeta \in \overline{\O}$
then we set
$$
\rho(z, \zeta) = |1 - z \cdot \overline{\zeta}|^{1/2} \, .
$$
\end{definition}

\begin{proposition} \sl
When $\Omega = B$, the unit ball, then the function $\rho$ is a metric on $\partial B$.  
For a more general smoothly bounded, strongly pseudoconvex domain, the function $\rho$
is a pseudometric.  That is to say, there is constant $C \geq 1$ such that
$$
\rho(z, \zeta) \leq C \bigl ( \rho(z, \xi) + \rho(\xi, \zeta) \bigr ) \, .
$$
\end{proposition}
{\bf Proof:}  The first assertion is Proposition 6.5.1 in [KRA5].  The second assertion is
proved in pp.\ 357--8 in [KRA1].
\endpf
\smallskip \\

\begin{proposition} \sl
The balls 
$$
\beta_2(z, r) = \{\zeta \in \O: \rho(z, \zeta) < r\} \, ,
$$
together with ordinary Euclidean volume measure $dV$, form a space
of homogeneous type in the sense of Coifman and Weiss [COW].
\end{proposition}
{\bf Proof:}  This is almost immediate from the preceding proposition,
but details may be found in Section 8.6 of [KRA1].  
\endpf
\smallskip \\

\begin{definition} \rm
For $z \in \O$ and $f \in L^1_{\rm loc}(\O)$ we define
$$
{\cal M}f(z) = \sup_{r > 0} \frac{1}{V(\beta_2(z, r)} \int_{\beta_2(z,r)} |f(\zeta)| \, dV(\zeta) \, .
$$
\end{definition}

\begin{theorem} \sl
The operator ${\cal M}$ is of weak type $(1,1)$ and of strong type $(p,p)$, $1 < p \leq \infty$.
\end{theorem}
{\bf Proof:}  Again this is a standard consequence of the previous proposition in the
context of spaces of homogeneous type.  See [COW].
\endpf
\smallskip \\

\begin{theorem} \sl
Let $\O$ be the unit ball $B$ in $\CC^n$.  Let $f$ be a locally
integrable function on $\O$.  Then there is a constant $C > 0$ such that,
for $z \in \O$,
$$
|{\cal B} f(z)| \leq C \cdot {\cal M} f(z) \, .
$$
\end{theorem}
{\bf Proof:}   It is easy to see that $|1 - z \cdot \overline{\zeta}| \geq (1/2)(1 - |z|^2)$.  Therefore
we may perform these standard estimates:
\begin{align}
|\B f(z)| & =  \left | \int_\O \B(z, \zeta) f(\zeta) \, dV(\zeta) \right | \notag \\
       & \leq  \sum_{j=-1}^\infty \int_{2^j (1 - |z|^2) \leq |1 - z \cdot \overline{\zeta}| \leq 2^{j+1} (1 - |z|^2)}
			     \B(z, \zeta) |f(\zeta)| \, dV(\zeta) \notag \\
       & \leq  \sum_{j=-1}^\infty \int_{|1 - z \cdot \overline{\zeta}| \leq 2^{j+1} (1 - |z|^2)}
			    \frac{(1 - |z|^2)^{n+1}}{[2^j(1 - |z|^2)]^{2n+2}} \, dV(\zeta) \notag \\
       & \leq  C \cdot \sum_{j=-1}^\infty 2^{-j(n+1)} \cdot \left [ \frac{1}{(1 - |z|^2)^{n+1} 2^{(j+1)(n+1)}} \right ] \int_{|1 - z \cdot \overline{\zeta}| \leq 2^{j+1} (1 - |z|^2)}
			 |f(\zeta)| \, dV(\zeta)  \notag \\
       & \leq  C \cdot \sum_{j=-1}^\infty 2^{-j(n+1)} \cdot \left [ \frac{1}{V(\beta_2(z, \sqrt{2^{j+1}(1 - |z|^2)})} \right ]
		\int_{\beta_2(z, \sqrt{2^{j+1}(1 - |z|^2)}}
			 |f(\zeta)| \, dV(\zeta)  \notag \\
\end{align}
The last line is majorized by
\begin{align}
       & \leq  C' \cdot \sum_{j=-1}^\infty 2^{-j(n+1)} {\cal M} f(z) \notag \\ 
       & \leq  C \cdot {\cal M} f(z) \, .  \tag*{$\BoxOpTwo$}  \\ \notag
\end{align}

\begin{theorem} \sl
Let $\O$ be the unit ball $B$ in $\CC^n$.  Let $f$ be an $L^p(\O, dV)$ function,
$1 \leq p \leq \infty$.
Then $\B f$ has radial boundary limits almost everywhere on $\partial \O$.
\end{theorem}
{\bf Proof:}  The proof follows standard lines, using Theorems 4.6 and 4.7.  See
the detailed argument in [KRA1, Theorem 8.6.11].
\endpf
\smallskip \\

In fact a slight emendation of the arguments just presented allow a more refined result.  

\begin{definition} \rm
Let $P \in \partial B$ and $\alpha > 1$.  Define the {\it admissible approach region
of aperture $\alpha$} by
$$
{\cal A}_\alpha(P) = \{z \in B: |1 - z \cdot \overline{\zeta}| < \alpha (1 - |z|^2) \} \, .
$$
\end{definition}

Admissible approach regions are a new type of region for Fatou-type theorems.  These
were first introduced in [KOR1], [KOR2] and generalized and developed in [STE] and
later in [KRA7].  Now we have

\begin{theorem} \sl
Let $f$ be an $L^p(B)$ function, $1 \leq p \leq \infty$.  Then, for
almost every $P \in \partial B$, 
$$
\lim_{{\cal A}_\alpha(P) \ni z \ra P} \B f(z) 
$$
exists.
\end{theorem}

In fact, using the Fefferman asymptotic expansion (as discussed in detail in the
next section), we may imitate the development of Theorems 4.6 and 4.7 and prove
a result analogous to Theorem 4.8 on any smoothly bounded, strongly pseudoconvex
domain.  We omit the details, as they would repeat ideas that we present elsewhere
in the present paper for slightly different purposes.

\section{Results on the Invariant Laplacian}
						     
If $g = (g_{jk})$ is a Riemannian metric on a domain $\Omega$
in complex Euclidean space, then there is a second-order partial
differential operator, known as the {\it Laplace-Beltrami
operator}, that is invariant under isometries of the metric.
In fact, if $g$ denotes the determinant of the metric matrix $g$, and if
$(g^{jk})$ denotes the inverse matrix, then
this partial differential operator is defined to be
$$
{\cal L} = \frac{2}{g} \sum_{j,k} \left \{ \frac{\partial}{\partial \bar{z}_j} \left ( g g^{jk} \frac{\partial}{\partial z_k} \right )
                + \frac{\partial}{\partial z_k} \left ( g g^{jk} \frac{\partial}{\partial \bar{z}_k} \right ) \right \} \, .
$$

Now of course we are interested in artifacts of the Bergman theory.  If $\Omega \ss \CC^n$ is a bounded
domain and $K= K_\Omega$ its Bergman kernel, then it is well known (see [KRA1]) that
$K(z,z) > 0$ for all $z \in \Omega$.  Then it makes sense to define
$$
g_{jk}(z) = \frac{\partial^2}{\partial z_k \partial \zbar_k} \log K(z,z)
$$
for $j, k = 1, \dots, n$.
Then Proposition 3.2 can be used to demonstrate that this metric---which is in fact a K\"{a}hler metric on $\Omega$---is 
invariant under biholomorphic mappings of $\Omega$.  In other words, any biholomorphic
$\Phi: \Omega \ra \Omega$ is an isometry in the metric $g$.  This is the celebrated {\it Bergman metric}.

If $\Omega \ss \CC^n$ is the unit ball $B$, then the Bergman kernel is given
by
$$
K_B(z, \zeta) = \frac{1}{V(B)} \cdot \frac{1}{(1 - z\cdot \zetabar)^{n+1}} \, ,
$$
where $V(B)$ denotes the Euclidean volume of the domain $B$.
Then
$$
\log K(z,z) = - \log V(B) - (n+1) \log (1 - |z|^2) .
$$
Further,
$$
\frac{\partial}{\partial z_j} \bigl ( - (n+1) \log (1 - |z|^2) \bigr )  = 
                                         (n+1) \frac{\bar{z}_j}{1 - |z|^2} 
$$
and
\begin{eqnarray*}
 \frac{\partial^2}{\partial z_j \partial \bar{z}_k} \bigl ( - (n+1) \log (1 - |z|^2 \bigr ) &  = &
                                                  (n+1) \left [ \frac{\delta_{jk}}{1 - |z|^2} + \frac{\bar{z}_j z_k}{(1 - |z|^2)^2} \right ] \\
              & = & \frac{(n+1)}{(1 - |z|^2)^2} \bigl [ \delta_{jk} (1 - |z|^2) + \bar{z}_j {z}_k \bigr ] \\
           & \equiv & g_{jk}(z) .
\end{eqnarray*}

When $n=2$  we have
$$
g_{jk}(z) = \frac{3}{(1 - |z|^2)^2} \bigl [ \delta_{jk}(1 - |z|^2) + \bar{z}_j z_k \bigr ] .
$$
Thus
$$
  \bigl ( g_{jk}(z) \bigr ) = \frac{3}{(1 - |z|^2)^2}
             \left ( \begin{array}{lr}
                 1 - |z_2|^2 & \bar{z}_1 z_2 \\
                 \bar{z}_2 z_1 & 1 - |z_1|^2
                     \end{array}
             \right ) .
$$
Let 
$$
\biggl ( g^{jk}(z) \biggr )_{j,k=1}^2
$$
represents the inverse of the matrix 
$$
\biggl ( g_{jk}(z) \biggr )_{j,k=1}^2 \ \ .
$$
Then an elementary computation shows that
$$
\biggl ( g^{jk}(z) \biggr )_{j,k=1}^2
     = \frac{1 - |z|^2}{3} \left ( \begin{array}{lr}
                         1  - |z_1|^2 & - z_2 \bar{z}_1 \\
                        - z_1 \bar{z}_2 & 1 - |z_2|^2 
                                   \end{array}
                           \right )
     = \frac{1 - |z|^2}{3} \bigl ( \delta_{jk} - \bar{z}_j z_k \bigr )_{j,k} .
$$
Let
$$
g \equiv \det \biggl ( g_{jk}(z) \biggr ) .
$$
Then
$$
 g = \frac{9}{(1 - |z|^2)^3} . 
$$

Now let us calculate.
If $\bigl ( g_{jk}\bigr )^2_{j,k=1}$ is the Bergman metric on the
ball in $\CC^2$ then we have
$$
\sum_{j,k} \frac{\partial}{\partial \bar{z}_j} \bigl ( g g^{jk} \bigr ) = 0
$$
and 
$$ 
\sum_{j,k} \frac{\partial}{\partial z_j} \bigl ( g g^{jk} \bigr ) = 0 .
$$

We verify these assertions in detail in dimension $2:$
Now
\begin{eqnarray*}
 g g^{jk} & = & \frac{9}{(1 - |z|^2)^3} \cdot \frac{1 - |z|^2}{3} (\delta_{jk} - \bar{z}_j z_k ) \\
          & = & \frac{3}{(1- |z|^2)^2} (\delta_{jk} - \bar{z}_j z_k ) .
\end{eqnarray*}
It follows that
$$
\frac{\partial}{\partial \bar{z}_j} \biggl [ g g^{jk} \biggr ] =
              \frac{6 z_j}{(1 - |z|^2)^3} \bigl ( \delta_{jk} - \bar{z}_j z_k \bigr )      
                                    - \frac{3z_k}{(1 - |z|^2)^2} .
$$
Therefore
\begin{eqnarray*}
          \sum_{j,k = 1}^2 \frac{\partial}{\partial \bar{z}_j} \biggl [ g g^{jk} \biggr ] & = & 
                   \sum_{j,k = 1}^2 \left [ \frac{6 z_j (\delta_{jk} - \bar{z}_j z_k)}{(1 - |z|^2)^3}
                             - \frac{3 z_j}{(1 - |z|^2)^2} \right ] \\
               & = & 6 \sum_k \frac{z_k}{(1 - |z|^2)^3} -
                     6 \sum_{j,k} \frac{|z_j|^2 z_k}{(1 - |z|^2)^3} 
                     - 6 \sum_k \frac{z_k}{(1 - |z|^2)^2}  \\
               & = &  6 \sum_{j} \frac{z_j}{(1 - |z|^2)^2}
                       - 6 \sum_k \frac{z_k}{(1 - |z|^2)^2} \\
               & = & 0 .
\end{eqnarray*}
The other derivative is calculated similarly.

Our calculations show that, on the ball in $\CC^2,$
\begin{eqnarray*}
{\cal L} & \equiv & \frac{2}{g} \sum_{j,k} \left \{ \frac{\partial}{\partial \bar{z}_j}
                                      \left ( g g^{jk} \frac{\partial}{\partial z_k} \right ) +
                                      \frac{\partial}{\partial z_k}
                                      \left ( g g^{jk} \frac{\partial}{\partial \bar{z}_j} \right ) \right \} \\
                 & = & 4 \sum_{j,k} g^{jk} \frac{\partial}{\partial \bar{z}_j} \frac{\partial}{\partial z_k} \\
                 & = & 4 \sum_{j,k} \frac{1 - |z|^2}{3} \bigl ( \delta_{jk} - \bar{z}_j z_k \bigr )
                           \frac{\partial^2}{\partial z_k \partial \bar{z}_j} .
\end{eqnarray*}

Now the interesting fact for us is encapsulated in the following proposition:

\begin{proposition} \sl
The Poisson-Szeg\H{o} kernel on the ball $B$ solves the Dirichlet problem
for the invariant Laplacian ${\cal L}$.  That is to say, if $f$ is a continuous
function on $\partial B$ then the function
$$
u(z) = \left \{ \begin{array}{lcr}
             \int_{\partial B} {\cal P}(z, \zeta) \cdot f(\zeta) \, d\sigma(\zeta) & \hbox{if} & z \in B \\  [.15in]
	       f(z)                                                                & \hbox{if} & z \in \partial B 
		\end{array}
       \right.
$$
is continuous on $\overline{B}$ and is annihilated by ${\cal L}$ on $B$.
\end{proposition}

This fact is of more than passing interest.  In one complex variable, the study of holomorphic
functions on the disc and the study of harmonic functions on the disc are inextricably linked
because the real part of a holomorphic function is harmonic and conversely.  Such is not
the case in several complex variables.  Certainly the real part of a holomorphic function
is harmonic.  But in fact it is more:  such a function is {\it pluriharmonic}.  For the
converse direction, any real-valued pluriharmonic function is locally the real part of a holomorphic
function.  This assertion is false if ``pluriharmonic'' is replaced by ``harmonic''.  

And the result of Proposition 5.1 should not really be considered to be surprising.  For the
invariant Laplacian is invariant under isometries of the Bergman metric, hence invariant under
automorphisms of the ball.  And the Poisson-Szeg\H{o} kernels behaves nicely under automorphisms.
E. M. Stein was able to take advantage of these invariance properties to give a proof of Proposition 5.1
using Godement's theorem---that any function that satisfies a suitable mean-value property
must be harmonic (i.e., annihilated by the relevant Laplace operator).  See [STE] for the details.
\smallskip \\

\noindent {\bf Sketch of the Proof of Proposition 5.1} \ \ 
Now
$$
{\cal L} u = {\cal L} \int_{\partial B} {\cal P}(z, \zeta) \cdot f(\zeta) \, d\sigma(\zeta)
            = \int_{\partial B} \biggl [ {\cal L}_z {\cal P} (z,\zeta)\biggr  ] \cdot f(\zeta) \, d\sigma(\zeta) \, .
$$
Thus it behooves us to calculate ${\cal L}_z {\cal P}(z, \zeta)$.   Now we shall calculate this
quantity for each fixed $\zeta$.  Thus, without loss of generality, we may compose with a unitary
rotation and suppose that $\zeta = (1 + i0, 0 + i0)$ so that (in complex dimension 2)
$$								   
\P = c_2 \cdot \frac{(1 - |z|^2)^2}{|1 - z_1|^4} \, .
$$
This will make our calculations considerably easier.
						   
By brute force, we find that

\small

\begin{align}
\frac{\partial \P}{\partial \zbar_1} & = -2(1 - z_1)(1 - |z|^2) \cdot \left [ \frac{-1 + z_1 + |z_2|^2}{|1 - z_1|^6} \right ] \notag \\ 
\frac{\partial^2 \P}{\partial \zbar_1 \partial z_1} & = \frac{-2}{|1 - z_1|^6} \cdot \left [ - |z_1|^2 - |z_1|^2|z_2|^2 + 3|z_2|^2 - z_1 |z_2|^2 \right. \notag \\ 
             & \null \ \ \ \ \quad \left. - 2|z_2|^4 - 1 + z_1 + \zbar_1 - \zbar_1|z_2|^2 \right ] \notag \\  
\frac{\partial^2 \P}{\partial \zbar_1 \partial z_2} & = \frac{-2(1 - z_1)}{|1 - z_1|^6} \cdot \left [ 2\zbar_2 
        - \zbar_2 z_1 - 2 \zbar_2 |z_2|^2 - \zbar_2 |z_1|^2 \right ] \notag \\  
\frac{\partial^2 \P}{\partial z_1 \partial \zbar_2} & = \frac{-2(1 - \zbar_1)}{|1 - z_1|^6} \cdot \left [ 2z_2 
        - z_2 \zbar_1 - 2 z_2 |z_2|^2 - z_2 |z_1|^2 \right ] \notag \\	 
\frac{\partial \P}{\partial z_2} & = \frac{-2z_2 + 2|z_1|^2 z_2 + 2|z_2|^2 z_2}{|1 - z_1|^4} \notag \\ 
\frac{\partial^2 \P}{\partial z_2 \partial \zbar_2} & = \frac{-2 + 2|z_1|^2 + 4|z_2|^2}{|1 - z_1|^4} \tag{5.1.1} \\  \notag
\end{align}

\normalsize

Now we know that, in complex dimension two,

\small

\begin{eqnarray*}
{\cal L}_z \P(z, \zeta) & = & \frac{4}{3}(1 - |z|^2) \cdot (1 - |z_1|^2) \cdot \frac{\partial^2 \P_z}{\partial z_1 \partial \zbar_1} +
	   \frac{4}{3}(1 - |z|^2) \cdot (- \zbar_1 z_2) \cdot \frac{\partial^2 \P_z}{\partial z_2 \partial \zbar_1} \\
	   && \quad + \frac{4}{3}(1 - |z|^2) \cdot (- \zbar_2 z_1) \cdot \frac{\partial^2 \P_z}{\partial z_1 \partial \zbar_2} +
	   \frac{4}{3}(1 - |z|^2) \cdot (1 - |z_2|^2) \cdot \frac{\partial^2 \P_z}{\partial z_2 \partial \zbar_2} \, .
\end{eqnarray*}

\normalsize

\noindent Plugging the values from (5.1.1) into this last equation gives

\small

\begin{eqnarray*}
{\cal L}_z \P(z, \zeta) & = & \frac{4}{3}(1 - |z|^2) \cdot (1 - |z_1|^2) \cdot 
	 \frac{-2}{|1 - z_1|^6} \cdot \biggl [ - |z_1|^2 - |z_1|^2|z_2|^2 \\
	  && \null \ \ \ \ \quad + 3|z_2|^2 - z_1 |z_2|^2 
              - 2|z_2|^4 - 1 + z_1 + \zbar_1 - \zbar_1|z_2|^2 \biggr ]  \\  
	  && \ + \frac{4}{3}(1 - |z|^2) \cdot (- \zbar_1 z_2) \\ 
	  && \quad \ \ \ \times  \frac{-2(1 - z_1)}{|1 - z_1|^6} \cdot \biggl [ 2\zbar_2 
        - \zbar_2 z_1 - 2 \zbar_2 |z_2|^2 - \zbar_2 |z_1|^2 \biggr ]  \\  
	  && \ + \frac{4}{3}(1 - |z|^2) \cdot (- \zbar_2 z_1) \\ 
	  && \quad \ \ \ \times  \frac{-2(1 - \zbar_1)}{|1 - z_1|^6} \cdot \biggl [ 2z_2 
        - z_2 \zbar_1 - 2 z_2 |z_2|^2 - z_2 |z_1|^2 \biggr ]  \\
	  && \ + \frac{4}{3}(1 - |z|^2) \cdot (1 - |z_2|^2) \cdot |1 - z_1|^2 \cdot
	    \frac{-2 + 2|z_1|^2 + 4|z_2|^2}{|1 - z_1|^6} \, . \\
\end{eqnarray*}

\def\z1{z_1}
\def\z1bar{\overline{z}_1}
\def\z2{z_2}
\def\z2bar{\overline{z}_2}

\normalsize

Multiplying out the terms, we find that

\small

\begin{eqnarray*}
{\cal L}_z \P(z, \zeta) & = & \frac{-2}{|1 - z_1|^6} \cdot \biggl [ - |z_1|^2  - 4|z_1|^2|z_2|^2 + 3|z_2|^2 - z_1 |z_2|^2 - 2|z_2|^4 - 1 \\
			&& \quad + z_1 + \zbar_1 - \zbar_1|z_2|^2 + |z_1|^4 + |z_1|^4|z_2|^2 + z_1 |z_1|^2|z_2|^2 \\
			&& \quad + 2|z_1|^2 |z_2|^4 + |z_1|^2 - z_1|z_1|^2 - \zbar_1|z_1|^2 + \zbar_1|z_1|^2|z_2|^2 \biggr ] \\
			&& \quad - \frac{2}{|1 - z_1|^6} \cdot \biggl [ -2 \zbar_1 |z_2|^2 + 3|z_1|^2 |z_2|^2 + 2|z_2|^4 \zbar_1 + \zbar_1 |z_2|^2 |z_1|^2 \\
			&& \quad - z_1 |z_1|^2 |z_2|^2 - 2 |z_1|^2 |z_2|^4 - |z_2|^2 |z_1|^4 \biggr ] \\
			&& \quad - \frac{2}{|1 - z_1|^6} \cdot \biggl [ -2 z_1 |z_2|^2 + 3|z_1|^2 |z_2|^2 + 2|z_2|^4 z_1 + z_1 |z_2|^2 |z_1|^2 \\
			&& \quad - \zbar_1 |z_1|^2 |z_2|^2 - 2 |z_1|^2 |z_2|^4 - |z_2|^2 |z_1|^4 \biggr ] \\
\end{eqnarray*}

\begin{eqnarray*}
			&& \quad - \frac{2}{|1 - z_1|^6} \cdot \biggl [ 1 - |z_1|^2 - 3|z_2|^2 + |z_1|^2 |z_2|^2 + 2|z_2|^4 - z_1 + z_1|z_1|^2 \\
			&& \qquad + 3 z_1 |z_2|^2 - z_1 |z_1|^2 |z_2|^2 - 2 z_1 |z_2|^4 - \zbar_1 + \zbar_1 |z_1|^2 + 3 \zbar_1 |z_2|^2 - \zbar_1 |z_1|^2 |z_2|^2 \\
			&& \qquad - 2 \zbar_1 |z_2|^4 + |z_1|^2 - |z_1|^4 - 3|z_1|^2 |z_2|^2 + |z_1|^4 |z_2|^2 + 2|z_1|^2 |z_2|^4 \biggr ] \, .	   \\
\end{eqnarray*}

\normalsize

And now if we combine all the terms in brackets a small miracle happens:  everything cancels.  The result is
$$
\L_z \P(z, \zeta) \equiv 0 \, .  \eqno \BoxOpTwo
$$
\vspace*{.15in}

Thus, in some respects, it is inappropriate to study holomorphic functions on the ball in $\CC^n$ using
the Poisson kernel.  The classical Poisson integral does {\it not} create pluriharmonic
functions, and it does not create functions that are annihilated by the invariant Laplacian.
In view of Proposition 5.1, the Poisson-Szeg\H{o} kernel is much more
apposite.  As an instance, Adam Koranyi [KOR] made decisive use of this observation
in his study (proving boundary limits of $H^2$ functions through admissible
approach regions ${\cal A}_\alpha$) of the boundary behavior of $H^2(B)$ functions.

It is known that the property described in Proposition 5.1 is special to
the ball---it is simply untrue on any other domain (see [GRA1], [GRA2] for
more detail on this matter). Now one of the points that we want to make in
this section is that the result of the proposition can be extended---in an
approximate sense---to a broader class of domains.

\begin{proposition} \sl
Let $\O \ss \CC^n$ be a smoothly bounded, strongly pseudoconvex domain and $\P$ its Poisson-\Sz\ kernel.
Then, if $f \in C(\partial \O)$ we may write
$$
\P f(z) = \P_1 f(z) + {\cal E} f(z) \, ,
$$
where
\begin{enumerate}
\item[{\bf (i)}]  The term $\P_1 f$ is ``approximately annihilated'' by the invariant Laplacian on $\Omega$;
\item[{\bf (ii)}]  The operator ${\cal E}$ is smoothing in the sense of pseudodifferential operators.
\end{enumerate}
\end{proposition}

\noindent We shall explain the meaning of {\bf (i)} and {\bf (ii)} in the course of the proofs of
these statements.
\smallskip \\

\noindent {\bf Proof of Proposition 5.2:}   We utilize of course the asymptotic
expansion for the \Sz\ kernel on a smoothly bounded, strongly pseudoconvex domain (see [FEF], [BMS]).
It says that, for $z, \zeta$ near a boundary point $P$, we have (in suitable biholomorphic
local coordinates)
$$
S_\Omega (z, \zeta) = \frac{c_n}{(1 - z \cdot \overline{\zeta})^n} + 
                 h(z,\zeta) \cdot \log  |1 - z \cdot \overline{\zeta}| \, .  \eqno (5.2.1)
$$
Here $h$ is a smooth function on $\overline{\O} \times \overline{\O}$.  

Now we calculate ${\cal P}(z, \zeta)$ in the usual fashion:
$$
\P_\Omega(z, \zeta) = \frac{|S(z, \zeta)|^2}{S(z, z)} = \frac{\displaystyle \biggl |\frac{c_n}{(1 - z \cdot \overline{\zeta})^n} +
                 h(z,\zeta) \cdot \log  |1 - z \cdot \overline{\zeta}| \biggr |^2}{\displaystyle \frac{c_n}{(1 - |z|^2)^n} 
                 + h(z,z) \cdot \log  (1 - |z|^2)} \, .  \eqno (5.2.2)
$$       
One can use just elementary algebra to simplify this last expression and obtain that, in suitable local coordinates near the boundary,
\begin{align}
\P_\Omega(z, \zeta) & =  c_n \cdot \frac{(1 - |z|^2)^n}{|1 - z \cdot \overline{\zeta}|^{2n}}  \notag	\\
   & \quad  \ + \ \frac{2(1 - |z|^2)^n}{|1 - z \cdot \overline{\zeta}|^n} \log |1 - z \cdot \overline{\zeta}| + 
	    {\cal O} \biggl [ (1 - |z|^2)^n \cdot \log |1 - z \cdot \overline{\zeta}| \biggr ]   \notag \\
   & \equiv   c_n \cdot \frac{(1 - |z|^2)^n}{|1 - z \cdot \overline{\zeta}|^{2n}} \ + \ {\cal E}(z,\zeta) \, .	 \tag{5.2.3} \\ \notag
\end{align}
     
Now the first expression on the righthand side of (5.2.3) is (in the local coordinates in which
we are working) the usual Poisson-\Sz\ kernel
for the unit ball in $\CC^n$.  The second is an error term which we now analyze.

In fact we claim that the error term is integrable in $\zeta$, uniformly in $z$, and the
same can be said for the gradient (in the $z$ variable) of the error term.  The first
of these statements is obvious, as both parts of the error term are clearly majorized
by the Poisson-Szeg\H{o} kernel itself.  As for the second part, we note that
the gradient of the error gives rise to three types of terms
\begin{align}
\nabla {\cal E} & \approx \frac{(1 - |z|^2)^{n-1}}{|1 - z \cdot \zetabar|^n} \cdot \log |1 - z \cdot \zetabar|  \notag \\  
                & \quad \ + \ \frac{(1 - |z|^2)^n}{|1 - z \cdot \zetabar|^{n+1}} \cdot \log |1 - z \cdot \zetabar| \notag \\  
		& \quad \ + \ \frac{(1 - |z|^2)^n}{|1 - z \cdot \zetabar|^{n+1}} \notag \\   
		& \ \equiv \ I + II + III \, .   \tag{5.2.4} \\  \notag
\end{align}
Now it is clear by inspection that $I$ and $II$ are majorized by the ordinary Poisson-Szeg\H{o} kernel, so
they are both integrable in $\zeta$ as claimed.  As for $III$, we must calculate:
\begin{eqnarray*}
\int_{\zeta \in \partial \Omega} \frac{(1 - |z|^2)^{n-1}}{|1 - z \cdot \overline{\zeta}|^{n+1}} \, d\sigma(\zeta) & \leq &
   \sum_{j=-1}^\infty \int_{2^j (1 - |z|^2) \leq |1 - z \cdot \overline{\zeta}| \leq 2^{j+1}(1 - |z|^2)}  \\
				  && \quad \frac{(1 - |z|^2)^{n-1}}{[2^j(1 - |z|^2)]^{n+1}} \, d\sigma(\zeta) \\
				  & \leq & \sum_{j=-1}^\infty \frac{1}{(1 - |z|^2)^2} \int_{|1 - z \cdot \overline{\zeta}| \leq 2^{j+1}(1 - |z|^2)} 
				   2^{-j(n+1)} \, d\sigma(\zeta) \\
				   & \leq &
   \sum_{j=-1}^\infty C \cdot \frac{2^{-j(n+1)}}{(1 - |z|^2)^2} \cdot \bigl [ \sqrt{2^{j+1}(1 - |z|^2)} \bigr ]^{2n-2} \\
                                  && \quad \times \bigl [ 2^{j+1} \cdot (1 - |z|^2) \bigr ] \\
				  & \leq & \sum_{j=-1}^\infty \frac{1}{(1 - |z|^2)^2} \cdot (1 - |z|^2)^{n-1} \cdot (1 - |z|^2) \\
				  && \quad \times 2^{-j(n+1)} \cdot 2^{(j+1)(n-1)} \cdot 2^{j+1} \\
				       & \leq & C \cdot 2^n (1 - |z|^2)^{n-2} \cdot \sum_{j=-1}^\infty  2^{-j} \\
				       & < & \infty \, .
\end{eqnarray*}

Thus we see that the Poisson-\Sz\ kernel for our strongly pseudoconvex domain $\O$ can be expressed,
in suitable local coordinates, as the Poisson-\Sz\ kernel for the ball plus an error
term whose gradient induces a bounded operator on $L^p$.  This means that the error term
itelf maps $L^p$ to a Sobolev space.  In other words, it is a smoothing operator (hence negligeble
from our point of view).

In fact there are several fairly well known results about the interaction
of the Poisson-Bergman kernel and the invariant Laplacian.  We summarize some
of the basic ones here.

\begin{proposition} \sl Let $f$ be a $C^2$ function on the unit
ball that is annihilated by the invariant Laplacian ${\cal
L}$.   Then, for any $0 < r < 1$ and $S$ the unit sphere,
$$
\int_S f(r \zeta) \, d\sigma(\zeta) = c(r) \cdot f(0) \, .
$$
Here $d\sigma$ is rotationally invariant measure on the sphere $S$.
\end{proposition}
{\bf Proof:}  Replacing $f$ with the average of $f$ over the orthogonal
group, this just becomes a calculation to determine the exact value
of the constant $c(r)$---see [RUD, p.\ 51].
\endpf
\smallskip 

\begin{proposition} \sl
Suppose that $f$ is a $C^2$ function on the unit ball $B$ that is annihilated
by the invariant Laplacian ${\cal L}$.  Then $f$ satisfies the identity
${\cal B}f = f$.  In other words, for any $z \in B$,
$$
f(z) = \int_B \B(z, \zeta) f(\zeta) \, dV(\zeta) \, .
$$
\end{proposition}
{\bf Proof:}  We have checked the result when $z = 0$ in the last
proposition.  For a general $z$, compose with a M\"{o}bius transformation
and use the biholomorphic invariance of the kernel and the differential
operator ${\cal L}$.
\endpf
\smallskip \\

\begin{remark} \rm
It is a curious fact (see [AFR]) that the converse of this last proposition
is only true (as stated) in complex dimensions $1, 2, \dots, 11$.  It is
false in dimensions 12 and higher.
\end{remark}			   

Finally we need to address the question of whether the
invariant Laplacian {\it for the domain $\Omega$} annihilates
the principal term of the righthand side of the formula
(5.2.3). The point is this. The biholomorphic change of
variable that makes (5.2.3) valid is {\it local}. It is valid
on a small, smoothly bounded subdomain $\O' \ss \O$ which
shares a piece of boundary with $\partial \O$. According to
Fefferman [FEF] (see also the work in [GRK1], [GRK2), there is
a smaller subdomain $\O'' \ss \O'$ (which also shares a piece
of boundary with $\partial \O$ and $\partial \O'$) so that the Bergman metric of
$\O'$ is close---in the $C^2$ topology---to the Bergman metric
of $\O$ {\it on the smaller domain $\O''$}. It follows then
that the Laplace-Beltrami operator $\L_{\O'}$ for the Bergman
metric of $\O'$ will be close to the Laplace-Beltrami operator
$\L_\O$ of $\O$ on the smaller subdomain $\O''$. Now, on
$\O'$, the operator $\L_{\O'}$ certainly annihilates the
principal term of (5.2.3). It follows then that, on $\O''$,
the operator $\L_\O$ {\it nearly} annihilates the principal
term of (5.2.3).  We shall not calculate the exact sense in
which this last statement is true, but leave details for the
interested reader.

This discussion completes the proof of Proposition 5.2.3.
\endpf
\smallskip \\

It is natural to wonder whether the Poisson-Bergman kernel
${\cal B}$ has any favorable properties with respect to important
partial differential operators.  We have the following positive
result:

\begin{proposition} \sl
Let $\O = B$, the unit ball in $\CC^n$, and $\B = {\cal B}_B(z,\zeta)$ its Poisson-Bergman
kernel.  Then ${\cal B}$ is plurisubharmonic in the $\zeta$ variable.
\end{proposition}
{\bf Proof:}  Fix a point $\zeta \in B$ and let $\Phi$ be an automorphism of $B$
such that $\Phi(\zeta) = 0$.  From Proposition 3.4 we then have
$$
\B_B(z, \zeta) = \B_B(\Phi(z), \Phi(\zeta)) \cdot |\hbox{det}\, J_\CC \Phi(\zeta)|^2 
   = \B_B(\Phi(z), 0) \cdot |\hbox{det}\, J_\CC \Phi(\zeta)|^2  \, .	\eqno (5.3.1)
$$
We see that the righthand side is an expression that is independent of $\zeta$ multiplied
times a plurisubharmonic function.   A formula similar to (5.3.1) appears in [HUA].

The same argument shows that $\B(\zeta, \zeta)$ is plurisubharmonic.
\endpf 
\smallskip \\

\section{Concluding Remarks}

The idea of reproducing kernels in harmonic analysis is an old one.
The Poisson and Cauchy kernels date back to the mid-nineteenth century.

Cauchy integral formula is special in that its kernel, which is
$$
\frac{1}{2\pi i} \cdot \frac{1}{\zeta - z} \, ,
$$
is just the same on any domain.  A similar statement is {\it not}
true for the Poisson kernel, although see [KRA7] for a study of
the asymptotics of this kernel.

The complex reproducing kernels that are indigenous to several complex variables
are much more subtle.   It was only in 1974 that C. Fefferman was able
to calculate Bergman kernel asymptotics on strongly pseudoconvex domains.
Prior to that, the very specific calculations of L. Hua [HUA] on concrete domains
with a great deal of symmetry was the standard in the subject.  A variant
of Fefferman's construction also applies to the Szeg\H{o} kernel (see also
[BMS]).  Carrying out an analogous program on a more general class of 
domains has proved to be challenging.

The current paper is an invitation to study yet another kernel---the Poisson-Bergman
kernel.  Inspired by the ideas of [HUA], this is a positive reproducing kernel
for the Bergman space.  There are many questions about the role of this new kernel
that remain unanswered.  We hope to investigate these matters in future work.

\newpage

\null \vspace*{.5in}

\section*{\sc References}
\vspace*{.2in}

\begin{enumerate}

\item[{\bf [AFR]}]  P. Ahern, M. Flores, and W. Rudin, An invariant volume-mean-value property,
{\it Jour.\ Functional Analysis} 11(1993), 380--397.

\item[{\bf [ARE]}]  J. Arazy and M. Engli\v{s}, Iterates and the boundary behavior of the Berezin transform, 
{\it Ann Inst.\ Fourier (Grenoble)} 51(2001), 1101--1133.

\item[{\bf [ARO]}] N. Aronszajn, Theory of reproducing kernels,
{\it Trans.\ Am.\ Math.\ Soc.} 68(1950), 337-404.

\item[{\bf [BEF]}] E. Bedford and J. E. Forn\ae ss, A
construction of peak functions on weakly pseudoconvex domains,
{\em Ann.\ of Math.} 107(1978), 555-568.

\item[{\bf [BER]}]  F. A. Berezin, Quantization in complex symmetric spaces,
{\it Math.\ USSR Izvestia} 9(1975), 341--379.

\item[{\bf [BCZ]}] C. A. Berger, L. A. Coburn, and K. H. Zhu,
Toeplitz operators and function theory in $n$-dimensions, {\it
Springer Lecture Notes in Math.} 1256(1987), 28--35.

\item[{\bf [BMS]}] L. Boutet de Monvel and J. Sj\"{o}strand,
Sur la singularit\'{e} des noyaux de Bergman et Szeg\"{o},
{\em Soc. Mat. de France Asterisque} 34-35(1976), 123-164.    

\item[{\bf [CAT]}] D. Catlin, Boundary behavior of holomorphic
functions on pseudoconvex domains, {\it J. Differential Geom.}
15(1980), 605--625.

\item[{\bf [COW]}] R. R. Coifman and G. Weiss, {\it Analyse
Harmonique Non-Commutative sur Certains Espaces Homogenes},
Springer Lecture Notes vol. 242, Springer Verlag, Berlin,
1971.

\item[{\bf [ENG1]}]  M. Engli\v{s}, Functions invariant under the Berezin
transform, {\it J. Funct. Anal.} 121(1994), 233--254.

\item[{\bf [ENG2]}]  M. Engli\v{s}, Asymptotics of the Berezin transform and 
quantization on planar domains, {\it Duke Math.\ J.} 79(1995), 57--76.

\item[{\bf [FEF]}] C. Fefferman, The Bergman kernel and biholomorphic
mappings of pseudoconvex domains, {\em Invent. Math.} 26(1974), 1-65.

\item[{\bf [GRA1]}] C. R. Graham, The Dirichlet problem for the Bergman
Laplacian. I, {\it Comm.\ Partial Differential Equations } 8(1983),
433--476.

\item[{\bf [GRA2]}] C. R. Graham, The Dirichlet problem for the Bergman
Laplacian. II, {\it Comm.\ Partial Differential Equations} 8(1983),
563--641.

\item[{\bf [GRK1]}]  R. E. Greene and S. G. Krantz, Stability
properties of the Bergman kernel and curvature properties of
bounded domains, {\it Recent Progress in Several Complex
Variables}, Princeton University Press, Princeton, 1982.

\item[{\bf [GRK2]}] R. E. Greene and S. G. Krantz, Deformation
of complex structures, estimates for the $\dbar$ equation, and
stability of the Bergman kernel, {\it Adv. Math.} 43(1982),
1-86.
		     
\item[{\bf [HUA]}] L. Hua, {\it Harmonic Analysis of Functions
of Several Complex Variables in the Classical Domains},
American Mathematical Society, Providence, RI, 1963.
	
\item[{\bf [KOR1]}]  A. Koranyi, Harmonic functions on Hermitian hyperbolic
space, {\it Trans.\ A. M. S.} 135(1969), 507-516.

\item[{\bf [KOR2]}]  A. Koranyi, Boundary behavior of Poisson integrals
on symmetric spaces, {\it Trans.\ A.M.S.} 140(1969), 393-409.

\item[{\bf [KRA1]}]  S. G. Krantz, {\it Function Theory of
Several Complex Variables}, $2^{\rm nd}$ ed., American
Mathematical Society, Providence, RI, 2001.

\item[{\bf [KRA2]}]  S. G. Krantz, {\it Explorations in Harmonic
Analysis with Applications to Complex Function Theory and the
Heisenberg Group}, Birkh\"{a}user Publishing, Boston, 2009, to
appear.

\item[{\bf [KRA3]}]  S. G. Krantz, {\it Cornerstones of
Geometric Function Theory: Explorations in Complex Analysis},
Birkh\"{a}user Publishing, Boston, 2006.

\item[{\bf [KRA4]}]  S. G. Krantz, A tale of three kernels,
{\it Complex Variables}, submitted.

\item[{\bf [KRA5]}]  S. G. Krantz, {\it Partial Differential
Equations and Complex Analysis}, CRC Press, Boca Raton, FL, 1992.

\item[{\bf [KRA6]}] S. G. Krantz, Calculation and estimation of the Poisson
kernel, {\it J. Math.\ Anal.\ Appl.} 302(2005), 143--148.

\item[{\bf [KRA7]}] S. G. Krantz, Invariant metrics and the boundary
behavior of holomorphic functions on domains in $\CC^{n}$, {\it Jour.\
Geometric.\ Anal.} 1(1991), 71-98.

\item[{\bf [PET]}]  J. Peetre, The Berezin transform and Ha-Plitz operators,
{\it J. Operator Theory} 24(1990), 165--186.

\item[{\bf [PHS]}] D. H. Phong and E. M. Stein, Hilbert
integrals, singular integrals, and Radon transforms. I, {\it
Acta Math.} 157(1986), 99--157.

\item[{\bf [RUD]}]  W. Rudin, {\it Function Theory in the Unit Ball of $\CC^n$},
Springer-Verlag, New York, 1980.

\item[{\bf [STE]}] E. M. Stein, {\it Boundary Behavior of
Holomorphic Functions of Several Complex Variables}, Princeton
University Press, Princeton, 1972.

\item[{\bf [ZHU]}]  K. Zhu, {\it Spaces of Holomorphic Functions in
the Unit Ball}, Springer, New York, 2005.

\end{enumerate}

\end{document}